\documentclass[12pt,a4paper]{article}
\usepackage{amssymb}
\usepackage{graphicx}
\newtheorem{dfn}{Definition}[section]
\newtheorem{tetel}[dfn]{Theorem}
\newtheorem{all}[dfn]{Proposition}
\newtheorem{lemma}[dfn]{Lemma}

\newcommand{\R}{\mathbb R}
\newcommand{\C}{\mathbb C}
\newcommand{\N}{\mathbb N}

\def\iM{{\cal M}}
\def\iP{{\cal P}}
\newcommand{\dint}{\int\kern-.05in\int}

\def\<{\langle}
\def\>{\rangle}

\newcommand{\ud}{\,d}
\newcommand{\ha}{\quad\textrm{if}\quad}

\newcommand{\tr}{\textrm{Tr\thinspace}}

\newcommand{\bv}{\nobreak\hfill $\square$}
\def\im{{\rm i}}
\begin{document}
\centerline{\LARGE{\bf Large deviation for the}}
\medskip
\centerline{\LARGE{\bf empirical eigenvalue density of}}
\medskip
\centerline{\LARGE{\bf truncated Haar unitary matrices}
\footnote{\today}}
\bigskip
\medskip\bigskip
\centerline{D\'enes Petz\footnote{E-mail: petz@math.bme.hu} 
    and J\'ulia R\'effy\footnote{E-mail: reffyj@math.bme.hu}}
\bigskip
\medskip
\centerline{Department for Mathematical Analysis}
\centerline{Budapest University of Technology and Economics}
\centerline{ H-1521 Budapest XI., Hungary}
\bigskip\bigskip

\begin{quote}
Let $U_m$ be an $m \times m$ Haar unitary matrix and $U_{[m,n]}$
be its $n \times n$ truncation. In this paper the large deviation
is proven for the empirical eigenvalue density of $U_{[m,n]}$
as $m/n \to \lambda $ and $n \to \infty$. The rate function and
the limit distribution are given explicitly. $U_{[m,n]}$ is the
random matrix model of $quq$, where $u$ is a Haar unitary in
a finite von Neumann algebra, $q$ is a certain projection and they
are free. The limit distribution coincides with the Brown measure
of the operator  $quq$. 
\end{quote}

\begin{quote}
\end{quote}

\begin{quote}
\end{quote}
\bigskip\bigskip
\noindent

\section{Introduction}

Although the asymptotics of the eigenvalue density of different random
matrices has been widely studied since the pioneering work of Wigner 
\cite{Wig}, the first large deviation theorem for the empirical 
eigenvalue density of self-adjoint Gaussian random matrices was proven 
by Ben Arous and Guionnet much later \cite{BAG}. After the publication 
of their work, several similar theorems were obtained for different kind 
of random matrices. In particular, Haar distributed unitaries were 
discussed by Hiai and Petz \cite{HPuni} and the monograph \cite{HP} 
contains more information about similar results (see also \cite{HPwig, BAZ}). 
Free probability theory has inspired non-commutative large deviation results
for random matrices recently, see \cite{Alice}, for example.
 
The aim of this article is to prove the large deviation theorem
for the empirical eigenvalue density of truncated Haar unitary 
random matrices, and to determine the limit measure. Let $U$ be an 
$m \times m$ Haar distributed unitary matrix. By truncating $m-n$
bottom rows and $m-n$ last columns, we get an $n\times n$ matrix. 
The truncated matrix is a contraction, hence the eigenvalues are in the 
unit disc. Our aim is to study the asymptotics of the empirical
eigenvalue density when $n \to \infty$ and $m/n \to \lambda$. The 
truncated Haar unitaries appeared in the works \cite{ZS, Collins}. 
Since our random matrix model is unitarily invariant, the limiting 
eigenvalue density is rotation invariant in the complex plane.
It turns out that the limiting density is supported
on the disc of radius ${1 / \sqrt{\lambda}}$. In the paper the large deviation
result is established and the exact form of the rate function is given.
The large deviation implies the weak convergence of the empirical
eigenvalue density of the truncated unitaries with probability one.

The paper is organized as follows. Section 2 contains some preliminaries
about potential theory and large deviations. The large deviation result is
stated in Section 3. Section 4 contains the proof of our main result.  In Section 5
of the paper we make a connection to free probability theory. The truncated 
Haar unitaries form a random matrix model for the non-commutative random 
variable $quq$, where $u$ is an appropriate unitary, $q$ is a projection 
and they are assumed to be free. We observe that the limiting eigenvalue 
density coincides with the Brown measure of the operator $quq$. Our
paper is based on the joint eigenvalue density of truncated unitaries.
In the Appendix we sketch the derivation of this formula following
the original paper \cite{ZS}.

\section{Preliminaries}

In this section we review the setting of large deviation for the
empirical eigenvalue density of random matrices and collect some useful
concepts and results from potential theory.

Assume that $T_n(\omega)$ is a random $n\times n$ matrix with complex
eigenvalues $\zeta_1(\omega), \dots , \zeta_n(\omega)$. 
(If we want, we can fix an ordering of the eigenvalues, for example, 
regarding their absolute values and phases, but that is not necessary.)
The {\it empirical eigenvalue density} of $T_n(\omega)$ is the random 
atomic measure
\[
P_n(\omega):=\frac{\delta(\zeta_1(\omega))+\dots
+\delta(\zeta_n(\omega))}{n}\, ,
\]
where $\delta(z)$ denotes the Dirac measure supported on $\{z\}\subset \C$. 
Therefore $P_n $ is a random measure, or a measure-valued random variable.

Let us recall the definition of the {\it large deviation principle} \cite{DZ}. 
Let $(P_n)$ be a sequence of measures on a topological space $X$. 
The large deviation principle holds with {\it rate function} $I: \iM(X) \to
\R^+ \cup \{+\infty\}$ in the scale $n^{-2}$ if
\[
\liminf_{n\to \infty} \frac{1}{n^2}\log P_n(G)\geq -\inf_{x\in G}I(x)
\]
for all open set $G\subset X$, and
\[
\limsup_{n\to \infty} \frac{1}{n^2}\log P_n(F)\leq -\inf_{x\in F}I(x)
\]
for all closed set $F\subset X$.

Let $U(m)$ be an $m \times m$ Haar distributed unitary matrix. 
By truncating $m-n$ bottom rows and $m-n$ last columns, we get a 
$n\times n$ matrix $U_{[m,n]}$. The truncated matrix $U_{[m,n]}$ is 
not a unitary but its operator norm is at most 1. Hence the
eigenvalues $\zeta_1 , \zeta_2 , \ldots , \zeta_n$ lie in
the disc $\mathcal D := \{ z \in \C \, : \, |z| \leq 1\}$. The
relevant topological space is $\mathcal M(\mathcal D)$, the space of
probability measures on  $\mathcal D$. Note that
this space is a compact metrizable space with respect to the weak
convergence of measures. Let $P_{[m, n]}$ be the
empirical eigenvalue density of $U_{[m,n]}$. Hence $P_{[m, n]}$
may be regarded as a measure on $\mathcal M(\mathcal D)$. 

We are going to benefit from the fact that the joint probability density
of the eigenvalues of $U_{[m,n]}$ is
$$
{1 \over C_{[m,n]}} \prod_{1\leq i < j\leq n} |\zeta_i -\zeta_j |^2 
\prod_{i=1}^n (1-|\zeta_i |^2)^{m-n-1}
$$
according to \cite{ZS}, see also the Appendix. The normalizing constant
\begin{equation}\label{E:konst}
C_{[m,n]}=\pi^nn! \prod_{j=0}^{n-1}{m-n+j-1\choose j}^{-1}
\frac{1}{m-n+j}
\end{equation}
was obtained in \cite{PDRJ}.

Next we recall some definitions and theorems of potential theory \cite{totik}.
For a signed measure $\nu$ on $\mathcal D$   
\[
\Sigma(\nu):=\dint_{\mathcal D^2}\log|z-w|\ud \nu(z)\ud\nu(w)
\]
is the negative logarithmic energy of $\nu$. Since
\[
\Sigma(\nu)=\inf_{\alpha<0}\dint_{\mathcal D^2}\max(\log|z-w|,\alpha)\ud \nu(z)
\ud\nu(w),
\]
this functional is upper semi-continuous. We want to show its concavity.

The following lemma is strongly related to the properties of the 
logarithmic kernel $K(z,w)=\log|z-w|$ (cf. Theorem 1.16 in \cite{landkof}).

\begin{lemma}\label{sigma}
Let $\nu$ be a compactly supported signed measure on $\C$ such that
$\nu(\C)=0$. Then $\Sigma(\nu)\leq 0$, and $\Sigma(\nu)=0$ if
and only if $\nu =0.$
\end{lemma}

From this lemma we can deduce strictly concavity of the functional $\Sigma$.
First we prove that 
\begin{equation}\label{E:midpoint}
\Sigma\left(\frac{\mu_1+\mu_2}{2}\right)\geq 
\frac{\Sigma(\mu_1)+\Sigma(\mu_2)}{2},
\end{equation}
for all $\mu_1,\mu_2\in\mathcal M(\mathcal D)$, moreover
the equality holds if and only if $\mu_1=\mu_2$. For this, apply
Lemma \ref{sigma} for the signed measure $\nu=\mu_1-\mu_2$. 
The strict midpoint concavity (\ref{E:midpoint}) implies strict concavity 
by well-known arguments.
 
Let $K \subset \C$ be a compact subset of the 
complex plane, and ${\mathcal M}(K)$ be the
collection of all probability measures with support in $K$.
The {\it logarithmic energy} $E(\mu)$ of a  $\mu\in{\mathcal M}(K)$
is defined as
\[
E(\mu):=\dint_{K^2}\log \frac{1}{|z-w|}\ud \mu(z)\ud\mu(w),
\]
and the {\it energy} $V$ of $K$ by
\[
V:=\inf\{E(\mu): \mu\in\mathcal M(K)\}.
\]
The quantity
\[
{\rm cap}(K):=e^{-V}
\]
is called the {\it logarithmic capacity} of $K$.
The {\it logarithmic potential} of $\mu\in\mathcal M(K)$ is the function 
\[
U^{\mu}:=\int_K\log \frac{1}{|z-w|}\ud \mu(w)
\]
defined on $K$. 
 
Let $K\subset \C$ be a closed set, and $Q:K\to (-\infty,\infty]$ be a
lower semi-continuous function which is finite on a set of positive
capacity. The integral  
\[I_Q(\mu):=\dint_{K^2}\log \frac{1}{|z-w|}\ud \mu(z)\ud \mu(w)+2\int_K
Q(z)\ud \mu(z)\]
is called {\it weighted energy}.

The following result  tells about the minimizer of the weighted 
potential (cf. Theorem I.3.3 in \cite{totik}).
 
\begin{all}\label{equi}
Let  $Q$ as above. Assume that $\sigma \in \mathcal M(K)$ has compact 
support, $E(\sigma)<\infty$ and 
\[
U^{\sigma}(z)+Q(z)
\]
coincides with a constant $F$ on the support of $\sigma$ and
is at least as large as $F$ on $K$. Then $\sigma$ is the  unique 
measure in $\mathcal M(K)$ such that
\[
I_Q(\sigma)=\inf_{\mu\in \mathcal M(K)}I_Q(\mu),
\]
i.e., $\sigma$ is the so-called {\it equilibrium  measure} associated with $Q$.
\end{all}

The following lemma is the specialization of Proposition \ref{equi} to  
a radially symmetric function $Q:\mathcal D\to (-\infty,\infty]$, i. e., 
$Q(z)=Q(|z|)$. We assume that $Q$ is differentiable on $(0,1)$ with
absolute continuous derivative bounded below,  moreover $rQ'(r)$
increasing on $(0,1)$ and 
 \[
\lim_{r\to 1}rQ'(r)=\infty.
\] 
Let $r_0\geq 0$ be the smallest number for which $Q'(r)>0$ for
all $r>r_0$, and we set $R_0$ be the smallest solution of $R_0Q'(R_0)=1$.
Clearly $0\leq r_0<R_0<1$.

\begin{lemma}\label{min}
If the above conditions hold, them the functional $I_Q$
attains its minimum at a unique measure $\mu_Q$ supported on the annulus
\[
S_Q=\{z:r_0\leq |z|\leq R_0\},
\]
and the density of $\mu_Q$ is given by
\[
\ud \mu_Q(z)=\frac{1}{2\pi}(rQ'(r))'\ud r\ud \varphi,\qquad z=re^{{\rm i}\varphi}.
\]  
\end{lemma}

\begin{biz}
The proof is similar to the one of Theorem IV.6.1 in \cite{totik}. 
Using the formula
\[
\frac{1}{2\pi}\int_0^{2\pi}\log \frac{1}{|z-re^{{\rm i}\varphi}|}\ud \varphi =
\left\{\begin{array}{ll}-\log r,&\textrm{if } |z|\leq r\\
-\log |z|,&\textrm{if } |z|> r,\end{array}\right.
\]
we get that   
\begin{eqnarray*}
U^{\mu}(z)&=&\frac{1}{2\pi}\int_{r_0}^{R_0}(rQ'(r))'
\int_0^{2\pi}\log\frac1{|z-re^{{\rm i}\varphi}|}\ud\varphi\ud r
\\
&=&Q(R_0)-\log R_0-Q(z),
\end{eqnarray*}
for  $z\in S_Q$, since $r_0=0$ or $Q'(r_0)=0$. We have 
\[
U^{\mu}(z)+Q(z)=Q(R_0)-\log R_0,
\]
which is clearly a constant.

Next we check that $U^{\mu}(z)+Q(z) \ge Q(R_0)-\log R_0$   for  $|z|< r_0$ and for
 $|z|>R_0$. So $\mu_Q$ satisfies conditions of Theorem \ref{equi} and it must
be the unique minimizer.
\end{biz}\bv

\section{The large deviation theorem}

Our large deviation theorem for truncated Haar unitaries is the
following.

\begin{tetel}\label{ldt} 
Let $U_{[m,n]}$ be the $n\times n$ truncation of an $m\times m$ 
Haar unitary random matrix and let  $1<\lambda<\infty$. If $m/n\to \lambda$ as 
$n\to \infty$, then the sequence of empirical eigenvalue densities
$P_n=P_{[m, n]}$ satisfies the large deviation principle in the scale
$1/n^{2}$ with rate function
\[
I(\mu):=-\dint_{\mathcal D^2}\log |z-w|\ud \mu(z)
\ud\mu(w)-(\lambda-1)\int_{\mathcal D}\log(1-|z|^2)
\ud \mu(z)+B,
\]
for $\mu \in \mathcal M(\mathcal D)$, where
\[
B:=-\frac{\lambda^2\log \lambda}{2}+\frac{\lambda^2\log (\lambda-1)}{2}-\frac{\log (\lambda-1)}{2}+
\frac{\lambda-1}{2}.
\]

Furthermore, there exists a 
unique $\mu_0\in \mathcal M(\mathcal D)$ given by the density  
\[
\ud\mu_0(z)=\frac{(\lambda-1)r}{\pi\left(1-r^2\right)^2}
\ud r\ud\varphi,\quad z=re^{\rm{i}\varphi}
\]
on $\{z:|z|\leq 1/\sqrt \lambda\}$ such that $I(\mu_0)=0$.
\end{tetel}

Set 
\[
F(z,w):=-\log |z-w|-\frac{\lambda-1}{2} \left(\log (1-|z|^2)
+\log(1-|w|^2)\right),
\]
and
\[
F_{\alpha}(z,w):=\min(F(z,w),\alpha),
\]
for $\alpha>0$. Since $F_{\alpha}(z,w)$ is bounded and continuous
\[
\mu\in \mathcal M(\mathcal D) \mapsto \dint_{\mathcal D^2}
F_{\alpha}(z,w) \ud \mu(z)\ud \mu(w).
\] 
is continuous in the weak* topology, when the support of $\mu$ is 
restricted to a compact set. The functional $I$ 
is written as
\begin{eqnarray*}
I(\mu)&&=\dint_{\mathcal D^2} F(z,w)\ud \mu(z)\ud \mu(w) +B\\
&&=\sup_{\alpha>0}
\dint_{\mathcal D^2} F_{\alpha}(z,w)\ud \mu(z)\ud \mu(w) +B\, ,
\end{eqnarray*}
hence $I$ is lower semi-continuous. 

\bigskip
We can write $I$ in the  form
\[
I(\mu)=-\Sigma(\mu)-(\lambda-1)\int_{\mathcal D}\log(1-|z|^2)
\ud \mu(z)+B.
\]
Here the first part $-\Sigma(\mu)$ is strictly convex (as it was established 
in the previous section) and the second part is affine in $\mu$. Therefore
$I$ is a strictly convex functional.

If $X$ is compact and $\mathcal A$ is a base for the topology, then the
large deviation principle is equivalent to the following
conditions (Theorem 4.1.18 in \cite{DZ}):
\[
-I(x)=\inf_{x\in G, G\in \mathcal A}
\left\{\limsup_{n\to \infty} \frac{1}{n^2} \log P_n(G)\right\}=
\inf_{x\in G, G\in \mathcal A}
\left\{\liminf_{n\to \infty} \frac{1}{n^2} \log P_n(G)\right\}
\]
for all $x\in X$.
We apply this result in the case $X=\mathcal M(\mathcal D)$, and we choose
$$
\left\{\mu'\in\mathcal M(\mathcal D):\left|\int_{\mathcal D}
z^{k_1}\overline z^{k_2}\ud \mu'(z)
-\int_{\mathcal D}z^{k_1}\overline z^{k_2}\ud \mu(z)
\right|<\varepsilon \textrm{ for } k_1+k_2\leq m\right\}.
$$
to be $G(\mu;m,\varepsilon)$. For $\mu \in\mathcal M(\mathcal D)$ 
the sets $G(\mu;m,\varepsilon)$ form a neighborhood base
of $\mu$ for the weak* topology of $\mathcal M(\mathcal D)$, 
where $m\in \N$ and $\varepsilon >0$.
To obtain the theorem, we have to prove that
\[
-I(\mu)\geq \inf_G\left\{\limsup_{n\to \infty} \frac{1}{n^2} \log P_n(G)\right\},\]
\[-I(\mu)\leq \inf_G\left\{\liminf_{n\to \infty} \frac{1}{n^2} \log
P_n(G)\right\},
\]
where $G$ runs over neighborhoods of $\mu$.

The large deviation theorem implies the almost sure weak convergence.

\begin{tetel}
Let $U_{[m,n]}$, $P_n$ and $\mu_0$ as in Theorem \ref{ldt}. Then 
\[
P_n(\omega)\stackrel{\scriptscriptstyle{n\to\infty}}{\longrightarrow}\mu_0\]
weakly with probability 1.
\end{tetel}

The proof is standard, one benefits from the compactness of the level
sets of the rate function and the Borel-Cantelli lemma is used, see
\cite{DZ}.

\section{Proof of the large deviation}

In this section we prove Theorem \ref{ldt}. Our method is based 
on the explicit form of the joint eigenvalue density. 

First we compute the limit of the normalizing  constant $C_{[m,n]}$
given in (\ref{E:konst}). 
\begin{eqnarray*}
B &=:&\lim_{n\to\infty}\frac{1}{n^2}\log C_{[m,n]}\\
&&=-\lim_{n\to\infty}\frac{1}{n^2}\sum_{j=1}^{n-1}\log {m-n+j-1\choose j}\\
&&=-\lim_{n\to \infty}\frac{1}{n-1}\sum_{i=1}^{n-1}\frac{n-1-i}{n-1}
\log\frac{m-n-1+i}{i}.
\end{eqnarray*}
Here the limit of a Riemannian sum can be recognized and this gives 
an integral:
\begin{eqnarray*}
B&=&-\int_0^1 (1-x)\log\left(\frac{\lambda-1+x}{x}\right)\ud x\\
&=&-\frac{\lambda^2\log \lambda}{2}+\frac{\lambda^2\log (\lambda-1)}{2}
-\frac{\log (\lambda-1)}{2}+\frac{\lambda-1}{2}.
\end{eqnarray*}

The lower and upper estimates are stated in the form of lemmas.

\begin{lemma}
For every $\mu \in \mathcal M(\mathcal D)$,
\[\inf_G\left\{\limsup_{n\to \infty} \frac 1{n^2}\log P_n(G)\right\}\leq
-\dint_{\mathcal D^2} F(z,w) \ud \mu(z)\ud \mu(w)-B\]
where $G$ runs over a neighborhood base of $\mu$.
\end{lemma}

This is the easier estimate, one can follow the proof of the earlier
large deviation theorems, see \cite{BAG, BAZ, HP}.
 
\begin{lemma}
For every $\mu \in\mathcal M(\mathcal D)$,
\[\inf_G\left\{\liminf_{n\to \infty} \frac 1{n^2}\log P_n(G)\right\}\geq
-\dint_{\mathcal D^2} F(z,w) \ud \mu(z)\ud \mu(w)-B,\]
where $G$ runs over a neighborhood base of $\mu$.
\end{lemma}
\begin{biz}
If
\[
\dint_{\mathcal D^2}F(z,w)\ud \mu(z)\ud \mu(w)
\]
is infinite, then we have a trivial case. Therefore we may 
assume that this double integral is finite.

Since $F(z,w)$ is bounded from below, we have
\[
\dint_{\mathcal D^2} F(z,w)\ud \mu(z)\ud \mu(w)=\lim_{k\to \infty}
\dint_{\mathcal D^2} F(z,w)\ud \mu_k(z)\ud \mu_k(w)
\]
with the conditional measure 
\[\mu_k(B)=\frac{\mu(B\cap \mathcal D_k)}{\mu(\mathcal D_k)},\]
for all Borel set $B$, where
\[
\mathcal D_k:=\left\{z:|z|\leq 1-\frac 1k\right\}.
\]
So it suffices to assume, that the support of $\mu$ is contained in
$\mathcal D_k$ for some $k\in \N$.

Next by possible regularization of the measure $\mu$, we we may assume 
that $\mu$ has a continuous density $f$ on the unit disc $\mathcal D$, 
and $\delta\leq f(z)$ for some $\delta>0$.

We want to partition the disc into annuli of equal measure. Let 
$k=[\sqrt n]$, and choose
\[
0=r_0^{(n)}\leq r_1^{(n)}\leq\dots \leq r_{k-1}^{(n)}\leq r_k^{(n)}=1,
\]
such that
\[
\mu\left(\left\{ z=re^{\im \varphi} : r\in [r_{i-1}^{(n)},r_i^{(n)}]\right\} \right)
=\frac{1}{k} \quad \textrm{for}\quad 1\leq i\leq k. 
\]
Note that \[k^2\leq n\leq k(k+2),\] and there exists a sequence $l_1,\dots, l_k$
such that 
$k\leq l_i\leq k+2$, for $1\leq i\leq k$, and $\sum_{i=1}^kl_i=n$. 
Now we partition radially. For fixed $i$ let
\[
0=\varphi_0^{(n)}\leq \varphi_1^{(n)}\leq\dots
\leq \varphi_{l_i-1}^{(n)}\leq \varphi_{l_i}^{(n)}=2\pi,
\]
such that
\[
\mu\left(\left\{z=re^{\im\varphi} : r\in [r_{i-1}^{(n)},r_i^{(n)}],\varphi\in 
[\varphi_{j-1}^{(n)},\varphi_j^{(n)}]\right\}\right)
=\frac{1}{kl_i} \quad \textrm{for}\quad 1\leq j\leq l_i.
\]
In this way we divided $\mathcal D$ into $n$ pieces, 
$S_1^{(n)},\dots , S_n^{(n)}$. Here
\begin{equation}\label{si}
\frac{\delta(1-\varepsilon_n)}{n}\leq\frac{\delta}{kl_i}=\int_{S_{i}^{(n)}}\ud z\leq
 \frac{1}{k^2\delta}\leq\frac{1+\varepsilon_n'}{n\delta},
 \end{equation}
where $\varepsilon_n=2/(\sqrt n+2)\to 0$ and $\varepsilon_n'=1/(\sqrt{n}-1)\to 0$ as $n\to \infty$. 
We can suppose, that
 \begin{equation}\label{atm}
 \lim_{n\to \infty}\left(\max_{1\leq i\leq n} {\rm{diam}} \left(S_{i}^{(n)}\right)\right)=0.
\end{equation}
In each part $S_{i}^{(n)}$ we take a smaller one $D_{i}^{(n)}$, similarly to 
$S_i^{(n)}$ by dividing the radial and phase intervals above into three 
equal parts, and selecting the middle ones, so that
\begin{equation}\label{di}
\frac{\delta(1-\varepsilon_n)}{9n}\leq\int_{D_{i}^{(n)}}\ud z
 \leq\frac{1+\varepsilon_n'}{9n\delta}.
 \end{equation}

We set
\[\Delta_n:=\left\{(\zeta_1,\dots,\zeta_{n}):\zeta_i\in D_{i}^{(n)},1\leq i\leq n\right\}.\]
For any neighborhood $G$ of $\mu$  
\[\Delta_n\subset \{\zeta \in \mathcal D^{n}:\mu_{\zeta}\in G\}\]
for every $n$ large enough. Then
\begin{eqnarray*}
\lefteqn{P_n(G)\geq \overline{\nu}_n(\Delta_n)}\\
&&=\frac{1}{Z_n} \int\dots \int_{\Delta_n} \exp \left(
(n-1)\sum_{i=1}^{n}(\lambda-1)\log\left(1-|\zeta_i|^2\right)\right)\\
&&\kern2.2in\times\prod_{1\leq i<j\leq n}|
\zeta_i-\zeta_j|^2\ud\zeta_1\dots \ud \zeta_{n}
\\
&&\geq \frac{1}{Z_n}\left( \frac{\delta(1-\varepsilon_n)}{9n}\right)
\exp\left((n-1)(\lambda-1)\sum_{i=1}^{n}\min_{\zeta \in  D_i^{(n)}}\log\left(1-|\zeta|^2
\right)\right)\\
&&\kern2.2in\times\prod_{1\leq i<j\leq n}\left(\min_{\zeta \in  D_i^{(n)},\eta \in D_j^{(n)}}
|\zeta-\eta|^2\right).
\end{eqnarray*}
Here for the first part we establish
\begin{eqnarray*}
\lim_{n\to \infty} &&\kern-.24in\frac{(n-1)(\lambda-1)}{n^2}\sum_{i=1}^{n} 
\min_{\zeta \in  D_i^{(n)}}\log\left(1-|\zeta|^2\right)\\
&&=\lim_{n\to\infty} \frac{\lambda-1}{n}\sum_{i=1}^n
\min_{\zeta \in  D_i^{(n)}}\log\left(1-|\zeta|^2\right)\\
&&=(\lambda-1)\int_{\mathcal D} \log\left(1-|\zeta|^2\right)f(\zeta)\ud \zeta,
\end{eqnarray*}
because of (\ref{atm}) and verify
\begin{eqnarray}\label{kettos}
\liminf_{n\to \infty} \frac{2}{n^2}\sum_{1\leq i<j\leq n} \log &&\kern-.3in 
\left( \min_{\zeta \in  D_i^{(n)},
\eta \in D_i^{(n)}}|\zeta-\eta|\right) \nonumber\\ & \geq &
\dint_{\mathcal D^2} f(\zeta)f(\eta)\log|\zeta-\eta|\ud \zeta\ud\eta.
\end{eqnarray}
for the second part.
\end{biz}\bv

The last step is to minimize $I$. Now we apply Lemma \ref{min} for
\[
Q(z):=-\frac{\lambda-1}{2}\log\left(1-|z|^2\right)
\]
on $\mathcal D$. This function satisfies the conditions of the lemma.
Hence the support of the limit measure $\mu_0$ is the disc
\[
S_{\lambda}=\left\{z:|z|\leq \frac{1}{\sqrt{\lambda}}\right\},
\] 
and the density is given by
\[
d\mu_0=\frac{1}{\pi}(rQ'(r))'\ud r\ud\varphi=\frac{1}{\pi}\frac{(\lambda-1)r}{\left(1-r^2\right)^2}
\ud r\ud\varphi,\quad z=re^{\rm{i}\varphi}.
\]

For this $\mu_0$ again by \cite{totik}
\begin{eqnarray*}
\lefteqn{I(\mu_0)=\frac 12Q\left(\frac{1}{\sqrt{\lambda}}\right)+\frac 1{2} \log \lambda
+\frac 1{2}\int_{S_{\lambda}} Q(z) d\mu_0(z)+B}\\
&&=-\frac{\lambda-1}{2}\log(\lambda-1)+\frac 1{2\lambda} \log \lambda-
\frac{(\lambda-1)^2}{2\pi}\int_0^{2\pi}
\int_0^{\frac{1}{\sqrt \lambda}} \frac{r\log(1-r^2)}{(1-r^2)^2}\ud r\ud \varphi\\
&&=-\frac{\lambda-1}{2}\log(\lambda-1)+\frac 1{2\lambda} \log \lambda-\frac {\lambda-1}{2}
\left( \lambda\log\left(\frac{\lambda-1}{\lambda}
\right)+1\right)+B=0.
\end{eqnarray*}
The uniqueness of $\mu_0$  satisfying $I(\mu_0)=0$ follows from the strict convexity of $I$.

\section{Some connection to free probability}

Let $Q_m$ be an $m\times m$ projection matrix of rank $n$, and let 
$U_m$ be an $m\times m$ Haar unitary. Then the matrix $Q_mU_mQ_m$ has 
the same non-zero eigenvalues as $U_{[m,n]}$, but it 
has $m-n$ zero eigenvalues. The large deviation result for $U_{[m,n]}$
is easily modified to have the following.

\begin{tetel}\label{ldt2} 
Let $1<\lambda<\infty$ and $Q_m,U_m$ as above. If $m/n\to \lambda$ as 
$n\to \infty$, then the sequence of empirical eigenvalue densities
$Q_m U_m Q_m$ satisfies the large deviation principle in the scale
$1/n^{2}$ with rate function
\[
\tilde I(\widetilde\mu):=
\left\{\begin{array}{lll}I(\mu),&\ha \widetilde\mu=(1-\lambda^{-1}) \delta_0+
\lambda^{-1} \mu,\\\phantom{p}\\
+\infty,&\textrm{ otherwise}
\end{array}\right.\]
Furthermore, the measure 
\[
\widetilde\mu_0= (1-\lambda^{-1})\delta_0+\lambda^{-1}\mu_0
\]
is the unique minimizer of $\tilde I$, and $\tilde I(\widetilde\mu_0)=0$.
\end{tetel}

\bigskip
Now let $\mathcal M$ be a von Neumann algebra and $\tau$ be a faithful
normal trace on $\mathcal M$. The pair $(\mathcal M, \tau)$ is often
called a {\it non-commutative probability space}. A unitary $u \in 
{\mathcal M}$
is called a {\it Haar unitary} if $\tau(u^k)=0$ for every non-zero integer
$k$. Let $q \in {\mathcal M}$ be a projection such that $\tau(q)=\lambda $. 
If $u$
and $q$ are {\it free} (see \cite{HP} or \cite{voi}
for more details about free probability),
then the above $(U_m, Q_m)$ is random matrix model of the pair $(u,q)$.
This means that 
$$
{1 \over m}E\left( \tr \iP (U_m, U_m^*, Q_m \right) 
\to \tau \left(\iP (u,u^*,q)\right)
$$
for any polynomial $\iP$ of three non-commuting indeterminants. This 
statement is a particular case of Voiculescu's fundamental result about
{\it asymptotic freeness} (\cite{Ve}, or Theorem 4.3.5 on p. 154 in \cite{HP}).

For an element $a$ of the von Neumann algebra $\iM$, the {\it Fuglede-Kadison}
determinant can be defined by:
$$
\Delta(a):=\lim_{\varepsilon \to +0}\exp \tau \left(\log(a^*a 
+\varepsilon I)^{1/2}\right).
$$
It was shown by L.G. Brown in 1983 that the function
$$
\lambda \mapsto \frac{1}{2\pi}\log \Delta (a-\lambda I)
$$
is subharmonic and its Laplacian (taken in the distribution sense) is
a probability measure $\mu_a$ concentrated on the
spectrum of $a$ \cite{B}. This measure is called the {\it Brown measure}
and it is a sort of extension of the spectral multiplicity measure of 
normal operators:
\begin{equation}
\tau(g(a)) = \int_{\C} g(z)\,d\mu_a(z)
\label{4}  
\end{equation}   
for any function $g$ on $\C$ that is analytic in a domain containing the 
spectrum of $a$. The Brown measure is computed for quite a few examples
in the paper \cite{Biane}. 

Let $u$ be a Haar unitary, and $q=q^*=q^2$ be free from $u$. Then $uq$
is a so-called {\it R-diagonal} operator and its Brown measure is
rotation invariant in the complex plane. According to \cite{haag} the Brown 
measure has an atom of mass $1-\lambda^{-1}$ at zero and the absolute 
continuous part has a density
\[
\frac{(\lambda-1)r}{\pi\lambda\left(1-r^2\right)^2}
\,d r\,d\varphi \qquad (z=re^{\rm{i}\varphi})
\]
on $\{z:|z|\leq 1/\sqrt \lambda\}$.
We just observe that this measure coincides with the limiting measure in
our large deviation theorem. In the moment we cannot deduce the Brown
measure from the large deviation result but it is definitely worthwhile to study
the relation.
 
\section*{Appendix}

Let $U_m$ be an $m \times m$ Haar unitary matrix and write it in the
block-matrix form
$$
\left( \begin{array}{cc}
A&B\\
C&D\end{array} \right), 
$$
where $A$ is an $n \times n$, $B$ is $n \times (m-n)$, $C$ is $(m-n)\times n$
and $D$ is an $(m-n) \times (m-n)$ matrix. The space of $n \times n$ (complex)
matrices is easily identified with $\R^{2n^2}$ and the push forward of the
usual Lebesgue measure is denoted by $\lambda_n$. It was obtained in \cite{Collins}
that for $m \ge 2n$, the distribution measure of the $n \times n$ matrix $A$
is absolute continuous with respect to $\lambda_n$ and the density is
\begin{equation}
C(n,m)\det(1-A^*A)^{m-2n}{\bf 1}_{\| A \| \le 1}d\lambda_n(A)\, .
\end{equation}

To determine the joint distribution of the eigenvalues $\zeta_1,\zeta_2,\dots, \zeta_n$
of $A$, we need only the matrices $A$ and $C$, and by a unitary transformation
we transform $A$ to an upper triangular form
\begin{equation}\label{E:AC}
\left( \begin{array}{ccccc}
\zeta_1 & \Delta_{1,2}& \Delta_{1,3}& \dots & \Delta_{1,n}\\
0 & \zeta_2 & \Delta_{2,3}& \dots & \Delta_{2,n} \\ . & . & .& \dots &.\\
0 & 0 & 0 & \dots &\zeta_n \\
C_1 & C_2 & C_3 & \dots & C_n
\end{array} \right),
\end{equation}
where $C_1,C_2,\dots,C_n$ are the column vectors of the matrix $C$.
It is well-known that the Jacobian of this  transformation is a
multiple of
$$
\prod_{1\le i < j \le n} |\zeta_i-\zeta_j|^2.
$$
 
Note that the columns of the matrix (\ref{E:AC}) are normalized and pairwise orthogonal.
Following the idea of  \cite{ZS}, we integrate out the variables $\Delta_{1,i},\Delta_{2,i}, \dots, 
\Delta_{i-1,i}, C_i$, $i \le n$.

One can construct  $(n-m)\times (n-m)$ matrices $X^{(i)}$ such that
\begin{equation}\label{E:faktor}
\Delta_{ij}=\frac{1}{\overline \zeta_i}C_i^*X^{(i)}C_j.
\end{equation}
We have  $X^{(1)}=I$ and
\[
X^{(i)}=I+\sum_{k<i}X^{(k)}\frac{C_k C_k^{\ast}}{|\zeta_k|^2}X^{(k)}.
\]
 
Since
\[
C_i^{\ast}C_i+\sum_{k<i}\overline\Delta_{ki}\Delta_{ki}=C_i^{\ast}X^{(i)}C_i,
\]
the vectors $C_i$ satisfy the equations
\begin{equation}\label{elli}
C_i^*X^{(i)}C_i=1-|\zeta_i|^2.
\end{equation}
Geometrically,  the point $(C_{1i},\dots, C_{m-n,i})$ lies in the ellipsoid given by $X^{(i)}$. 
To compute the volume of this ellipsoid it is enough to know the determinant of $X^{(i)}$
and this is obtained from the above recursion:
\[
\det X^{(i)}=\frac{\det X^{(i-1)}}{|\zeta_{i-1}|^2}=\prod_{j<i}\frac{1}{|\zeta_j|^2}.
\]

After this preparation we move to integration.  First we integrate with respect to the last column $
\Delta_{1,n},\Delta_{2,n}, \dots, \Delta_{n-1,n}C_n$.
For fixed $\Delta_{1,n}\dots \Delta_{n-1,n}$ the distribution of $C_{1,n},\dots, C_{{m}-n-1,n}$ is uniform on the set  
\[|C_{1,n}|^2+\dots +|C_{{m}-n-1,n}|^2\le 1-|\zeta_n|^2-|\Delta_{1,n}|^2\dots |\Delta_{n-1,n}|^2,\]
i.e. inside the ellipsoid defined by (\ref{elli}).
The volume of this ${m}-n-1$ dimensional complex ellipsoid is 
\[
\frac{(1-|\zeta_n|^2)^{{m}-n-1}}{\det X^{(n)}}=(1-|\zeta_n|^2)^{m-n-1}
\prod_{i<n}|\zeta_i|^2,
\]
Integration out of $\Delta_{i,n}$ gives a factor $|\zeta_i|^{-2}$ from (\ref{E:faktor}) 
and all together
we obtain a factor $(1-|\zeta_n|^2)^{m-n-1}$ from the last column.  
The same procedure may be applied to the other columns.

\section*{Acknowledgement}
The authors thank to Fumio Hiai, Ofer Zeitouni and the referee for useful remarks. The work has been
partially supported by Hungarian OTKA No. T032662 .

\end{document}